\documentclass{amsart}
\usepackage{amssymb,amsfonts,latexsym}
\usepackage{graphics,verbatim}
\usepackage{graphicx}
\usepackage{hyperref}
\newtheorem{theorem}{Theorem}[section]
\newtheorem{proposition}{Proposition}[section]
\newtheorem{lemma}{Lemma}[section]

\usepackage{color}

 \newcommand{\<}{\left\langle}
\renewcommand{\>}{\right\rangle}

\newcommand{\eps}{\varepsilon}

\newcommand{\To}{\longrightarrow}

\newcommand{\norm}[1]{\left\Vert#1\right\Vert}

\newcommand{\be} {\begin{equation}}
\newcommand{\ee} {\end{equation}}
\newcommand{\bea} {\begin{eqnarray}}
\newcommand{\eea} {\end{eqnarray}}
\newcommand{\Bea} {\begin{eqnarray*}}
\newcommand{\Eea} {\end{eqnarray*}}
\newcommand{\pa} {\partial}

\newcommand{\al} {\alpha}
\newcommand{\ba} {\beta}
\newcommand{\de} {\delta}

\newcommand{\ga} {\gamma}
\newcommand{\Ga} {\Gamma}
\newcommand{\Om} {\Omega}

\newcommand{\De} {\Delta}
\newcommand{\la} {\lambda}

\newcommand{\nequiv} {\not\equiv}

\newcommand{\var} {\varepsilon}
\newcommand{\f}{\frac}
\newcommand{\R}{\mathbb R}
\newcommand{\N}{\mathbb N}
\newcommand{\Rn}{\mathbb R^N}
\newcommand{\Iom}{\int_{\Omega}}
\newcommand{\deb}{\rightharpoonup}

\newcommand{\dx}{\,\mathrm{d}x} 
\newcommand{\dy}{\,\mathrm{d}y} 
\newcommand{\dxdy}{\, \mathrm{d}x \,\mathrm{d}y} 

\catcode`\@=11

\makeatletter \@addtoreset{equation}{section} \makeatother

\begin{document}

\title[On multiplicity results for fractional Brezis Nirenberg type equations]{On the
	existence of multiple solutions for fractional Brezis Nirenberg type equations}


\author{Debangana Mukherjee}
\address{Department of Mathematics and Statistics, Masaryk University, 61137 Brno,  Czech Republic}

\email{mukherjeed@math.muni.cz, \,\, debangana18@gmail.com }

\date{}
\maketitle

	\begin{abstract} The present paper studies the non-local fractional analogue of the famous  paper of Brezis and Nirenberg in \cite{Brezis-Nirenberg}. Namely, we focus on the following model,
		\begin{align*}
		\left(\mathcal{P}\right)
		\begin{cases}
		\left(-\Delta\right)^s u-\la u &= \al |u|^{p-2}u +  \ba|u|^{2^*-2}u \quad\mbox{in}\quad \Omega,\\
		u&=0\quad\mbox{in}\quad\mathbb{R}^N\setminus\Omega,
		\end{cases}
		\end{align*}
	where $(-\De)^s$ is the fractional Laplace operator, $s \in (0,1)$,
	with $N \geq 3s$, $2<p<2^*$, $\ba>0, \la, \al \in \R$ and establish the existence of nontrivial solutions and sign-changing solutions for the problem $(\mathcal{P})$.

\end{abstract}

\tableofcontents
\section{Introduction}
In the famous paper of Brezis and Nirenberg \cite{Brezis-Nirenberg}, they have researched on the following nonlinear critical elliptic partial differential equation:
\begin{align}\label{eq:BN}
\begin{cases}
-\Delta u-\la u = \al |u|^{p-2}u +  |u|^{2^*-2}u \quad\mbox{in}\quad \Omega,\\
u>0\quad\mbox{in}\quad \Omega, u=0\quad\mbox{on}\quad \partial \Omega,
\end{cases}
\end{align}
where $\Omega$ is a smooth bounded domain in $\Rn$, $2<p<2^*$ where $2^*=\frac{2N}{N-2}(N \geq 3)$ is the Sobolev critical exponent, $\la,\al$ are parameters. They have proved the following (see Corollaries 2.1-2.4 in \cite{Brezis-Nirenberg}):
\textit{
	\begin{itemize}
		\item [(I)]
		For $N \geq 4$,
		problem (\ref{eq:BN}) has a positive nontrivial solution for $\al>0$ and $\la \in (0,\la_1)$ where $\la_1$ is the first eigenvalue of $-\Delta$ with Dirichlet boundary condition in $\Omega$.
		\item [(II)] 
		For $N=3$, $\al>0$, $\la \in (0,\la_1)$, problem (\ref{eq:BN}) has a solution provided $4<p<6$ and problem (\ref{eq:BN}) has a solution for each $\al \geq \al_0$ for some $\al_0>0$ if $2<p \leq 4$.
	\end{itemize}	
}
In the interest of the pioneering work of H.R. Brezis and L. Nirenberg \cite{Brezis-Nirenberg}, a massive study is doing the rounds about the results dealing with semilinear problems involving critical Sobolev exponents, we mention some of the well-celebrated papers \cite{ Cerami-94, Fortunato-85, Fortunato-84, Struwe-86} and the references therein. In \cite{Bandle-02}, Bandle and Benguria have studied the classical Brezis-Nirenberg problem \cite{Brezis-Nirenberg} on $\mathbb{S}^3$.  In \cite{Miyagaki-19}, the authors have worked on the Brezis-Nirenberg problem on the hyperbolic space. In \cite{Benjin-07}, the authors have researched about the existence of non-trivial solutions to semilinear Brezis-Nirenberg problems involving Hardy potential and singular coefficients.

In the present paper, we are interested in the following nonlocal elliptic equation with Sobolev-critical exponent:
\begin{align*}
\left(\mathcal{P}\right)
\begin{cases}
\left(-\Delta\right)^s u-\la u &= \al |u|^{p-2}u +  \ba|u|^{2^*-2}u \quad\mbox{in}\quad \Omega,\\
u&=0\quad\mbox{in}\quad\mathbb{R}^N\setminus\Omega,
\end{cases}
\end{align*}
where $s\in(0,1)$ is fixed, $N \in \N$, $N>2s,\,2^*=\frac{2N}{N-2s}$ (the fractional critical Sobolev exponent), $2<q<2^*, \al, \ba>0, \lambda \in \R$, $\Om \subset \R^N$ is a bounded domain with smooth boundary
 and $(-\Delta)^s$ is the fractional Laplace operator, which (up to normalization factors) may be defined as:
\begin{align*}
-\left(-\Delta\right)^s u(x)=\frac{1}{2}\int_{\mathbb{R}^N}\frac{u(x+y)+u(x-y)-2u(x)}{|y|^{N+2s}}\dy,\,\,\,x\in\mathbb{R}^N.
\end{align*}

In the nonlocal setting, Brezis-Nirenberg type problems  are  in the pipeline and have been widely studied by many researchers. To name a few, we cite \cite{Chen-17}, \cite{Miyagaki-16}, \cite{Squassina-16},  \cite{Servadei-Valdinoci-13},\cite{Servadei-Valdinoci-15*}, \cite{Servadei-Valdinoci-15}. 
In \cite{Servadei-Valdinoci-15}, Servadei and Valdinoci have studied the following model:
\begin{align}\label{eq:Ser-Val}
\begin{cases}
(-\De)^s u-\la u = |u|^{2^*-2}u \quad\mbox{in}\quad \Omega,\\
 u=0\quad\mbox{in}\quad \Rn \setminus \Omega,
\end{cases}
\end{align}
where $\Om$ is an open bounded set of $\Rn$ with Lipschitz boundary, $N \geq 4s$, $s \in (0,1)$. They have proved that for any $\la \in (0, \la_{1,s})$, problem (\ref{eq:Ser-Val})
admits a  nontrivial solution $u \in H^s(\Rn)$ such that $u=0$ a.e. in $\Rn \setminus \Om$, where $\la_{1,s}$ is the first eigenvalue of the nonlocal operator $(-\De)^s$ with homogeneous Dirichlet boundary condition. 
In \cite{Tan-11}, the author has dealt with the existence and nonexistence of positive solutions to Brezis-Nirenberg type problems which involve the square root of the Laplace operator in a smooth bounded domain in $\Rn$.
In \cite{Miyagaki-16}, the authors have worked on 
the existence, nonexistence, and regularity of weak solutions for a non-local system involving fractional Laplacian. Very recently, in \cite{Colorado-19}, Colorado and Ortega have studied existence of solutions to the nonlocal critical Brezis-Nirenberg problem involving mixed Dirichlet-Neumann boundary conditions. In \cite{Cora-18}, Cora and Iacopetti have studied the asymptotic behavior and qualitative properties of least energy radial sign-changing solutions of the problem (\ref{eq:Ser-Val}) in a ball of $\Rn$.

\textbf{Functional Setting}. We mean by $H^s(\Rn)$ the usual fractional Sobolev space endowed with the so-called Gagliardo norm
$$\norm{g}_{H^s(\Rn)}=\norm{g}_{L^2(\Rn)}+\bigg(\int_{\R^{2N}} \frac{|g(x)-g(y)|^2}{|x-y|^{N+2s}}\dxdy\bigg)^{1/2}.$$ 
Let us signify $Q:=\R^{2N} \setminus (\Om^c \times \Om^c), \, \Om^c=\Rn \setminus \Om$ and
we define $$X_0 :=\Big\{u \in H^s(\mathbb{R}^N):u=0 \quad\text{a.e. in}\quad \Rn \setminus \Om\Big\}, $$ with the norm   
$$\norm{u}_{X_0}=\left(\int_{\mathbb{R}^{2N}}\frac{|u(x)-u(y)|^2}{|x-y|^{N+2s}}\dxdy\right)^{1/2}.$$
With this norm, $X_0$ is a Hilbert space with the scalar product
$$\<u,v\>_{X_0}=\int_{\mathbb{R}^{2N}}\frac{(u(x)-u(y))(v(x)-v(y))}{|x-y|^{N+2s}}\dxdy, $$
( see \cite[lemma 7]{Ser-Val-12}). For further details on $X_0$ and for their properties, we refer to \cite{NePaVal} and the references therein.  

Define $$\mathcal{D}^{s,2}(\Rn):=\{ u \in L^{2^*}(\Rn) : \int_{\R^{2N}} \frac{|u(x)-u(y)|^2}{|x-y|^{N+2s}} \dxdy < \infty  \},$$
and 
\begin{gather}\label{Eq:S}
{S}:=\inf_{u \in \mathcal{D}^{s,2}(\Rn) \setminus \{0\}}  \frac{\displaystyle\int_{\R^{2N}}{\frac{|u(x)-u(y)|^2}{|x-y|^{N+2s}}\dxdy }}{(\int_{\Rn}{ |u(x)|^{2^*}\dx})^{\frac{1}{2^*}}},
\end{gather}
is the best Sobolev constant for the fractional Sobolev embedding.

In this article, we will prove the following main results :

\begin{theorem} \label{thm.1}
	Under the assumptions $N \in \N$ with $N \geq 4s,\, s \in (0,1) \, \la,\al,\ba \in (0,\infty)$,  $\mathcal{(P)}$ has a non-trivial solution and has a sign changing solution provided $\la \geq \la_1$, where $\la_1$ is the first eigen value of $\left(-\Delta\right)^s$ in $X_0$.
\end{theorem}

\begin{theorem} \label{thm.2}
	Let $N \geq 3s, \la \in \R, \al,\ba>0. $ Then for all $m \in \N, $ there exists $\ba_m \geq 0$ such that $\mathcal{(P)}$ has $m$ nontrivial 
	solutions for all $\ba \in (0,\ba_m). $ Moreover, for $\la \geq \la_1, \mathcal{(P)}$ has $m$ sign changing solutions. 
\end{theorem}

\begin{theorem}\label{thm.3}
	Let $N \geq 4s$ with $s \in (0,1)$, $\la >0,\al<0,\ba>0$. Then, problem $(\mathcal{P})$ has a non-trivial solution for all $\al \in (-\al_0,0)$ for some $\al_0>0$.
\end{theorem}

In the classical case, this present problem is addressed in \cite{Yue-Zou}. The nonlocal framework has made this problem challenging. As far as we know, such result for existence of multiple and sign-changing solutions in the non-local framework, is not available in the literature.

\textbf{Organization of the paper}. The present manuscript consists of the following sections. Section 2 is devoted to recall  some preliminary results. In section 3, we prove Theorem \ref{thm.1}. Section 4 consists of the proof of Theorem \ref{thm.2}. In section 5, we establish the proof of Theorem \ref{thm.3}.

\section{Preliminaries}

\subsection{Eigenvalue problem}
In this section, we focus on the following eigenvalue problem:
\begin{align}\label{Eq:pblm-la}
\left(\mathcal{E}\right)
\begin{cases}
\left(-\Delta\right)^s u &= \la u \quad\mbox{in}\quad \Omega,\\
u&=0\quad\mbox{in}\quad\mathbb{R}^N\setminus \Omega,
\end{cases}
\end{align}
where $s \in (0,1), N>2s$ and $\Om \subset \Rn$ be an open, smooth bounded domain. We say $\la \in \R$ is an eigen value of $(-\De)$ if there exists a nontrivial solution of $(\mathcal{E})$ in $X_0.$ We recall the following properties of eigenvalues of $-\De$ in $X_0$. We refer [Proposition 9, \cite{Ser-Val-13}] for details.

\begin{itemize}
\item [(i)]	
The problem $(\mathcal{E})$ has an eigenvalue $\la_1>0$ which can be characterised as:
 $$ \la_1:=\inf_{u \in X_0 \setminus \{0\} } \frac{\|u\|_{X_0}}{|u|_2}, $$
where $|\cdot|_{p}$ denotes the $L^p$ norm in $\Om$ for $1 \leq p \leq \infty$. Moreover,	$\la_1$ is simple and the eigen function corresponding to $\la_1$ is not sign-changing.
\item [(ii)]
The set of eigen values of $(\mathcal{E})$ consist of a sequence $\{\la_k  \}_{k \geq 1}^{\infty}$ such that $0< \la_1<\la_2 \leq \la_3\leq \cdots \leq \la_k \leq \la_{k+1}\leq \cdots$ and $\lim_{k \to \infty}\la_k=\infty$.
\item [(iii)]
The sequence $\{e_k\}_{k \in \N}$ of eigen-functions in $X_0$ corresponding to $\{\la_k  \}_{k \in \N}$ forms an orthonormal basis of $L^2(\Om)$ and an orthogonal basis of $X_0$.
\item [(iv)]
For each $k \in \N$, the eigen value $\la_k$ has finite multiplicity.
\end{itemize}	
For $k \geq 1$,
let us denote $X_k:= \mbox{span}\{e_1, e_2, \cdots, e_k\}$. 
For $\la>0$, let us denote $$\la^+=\min_{k \in \N}   \{\la_k:\la_k>\la \}. $$ 
We note that
\begin{gather}\label{eq:lambda}
 \la^+=\la_l \, \text{for some} \, l \in \N \, \text{and}\, \la<\la_l. 
 \end{gather}
 For $\var > 0$, let us consider  the functions $v_\var :\Rn \to \R$ defined by:
 \begin{gather}\label{Eq:v-eps}
 v_\var(x)=\frac{[N(N-2s)\var]^{\frac{N-2s}{4}}}{(\var+|x|^2)^{\frac{N-2s}{2}}} \quad\text{for all}\quad x \in \Rn. 
 \end{gather}
 We know that (see \cite{CT-04})
 $S$, defined in (\ref{Eq:S}), is achieved by the family $\{v_\var\}_{\var>0}$. 
 Let us consider a smooth function $\psi$ of 
 $B_{1/2}(0)$ w.r.t $B_1(0), $ such that 
 	$\psi \in C_c^{\infty}(\Rn)$, $0 \leq \psi \leq 1$, $\psi(x) \equiv 1 \, \text{on}\, \text B_{1/2}(0)$ and $\mbox{supp}(\psi) \subset B_1(0)$.
 
Let us define $\phi_\var : \Rn \to \R$ by
$$\phi_\var(x)=\psi(x)v_\var(x), \, x \in \Rn. $$ 
 For $k \in \N$, let us denote the eigen space corresponding to the eigen value $\la_k$ by $Y_k$, that is,
 \begin{gather*}
 Y_k:= \big\{ u \in X_0 : (-\De)^su=\la_k  u \, \text{in}\, \Om  \big\}.
 \end{gather*}
 By property (iv) above, we note that $\text{dim}Y_k < \infty$ for each $k \in \N$.
Let us consider the orthogonal projection $P_k:X_0 \to Y_k$ of $X_0$ onto $Y_k$.
Set ${\eta_\var}:=(I-P_{l-1})\phi_\var=\phi_\var-P_{l-1}\phi_\var $  
where $l$ is given in (\ref{eq:lambda})
and define the set $V_\eps$ by
 \begin{gather}\label{eq:V-eps}
 V_\var:=\Big\{u \in X_0 :u=v+t {\psi_\var}, v \in X_{l-1},t \in \R \Big\} ,
 \end{gather}
  where
\begin{align}
{\psi_\var}=\begin{cases}
{\phi_\var}  \quad\text{if}\quad \la \neq \la_{l-1}, \\
{\eta_\var} \quad\text{if}\quad  \la =\la_{l-1}.
\end{cases}
\end{align}
We observe that $V_\var=X_{l-1}+\mbox{span} \{{\psi_\var}\}$.

We need the following result from  [Proposition 12, \cite{Servadei-Valdinoci-13}] to establish our main result Theorem \ref{thm.1}.
\begin{proposition}
For $\la \in (0,\infty)$, let us define the functional $Q_\lambda:H^s(\mathbb{R}^N)\setminus\{0\}\to \mathbb{R}$  by
\begin{gather*}
Q_\lambda(u):=\frac{\|{u}\|^2_{X_0}-\la|u|_{2}^2}{|u|_{2^*}^2}\quad\mbox{for all}\quad u\in H^s(\mathbb{R}^N) \setminus \{0\}.
\end{gather*}
Then, there exists $\eps_0>0$ (small enough) such that for all $\eps \in (0,\eps_0)$ and for all $\la \in (0,\infty)$, 
	\begin{gather}\label{Eq:Q-la}
	\sup_{u \in V_\var}Q_\lambda(u) < S,
	\end{gather}
where $S$ is defined in \eqref{Eq:S}.		
\end{proposition}

\section{Proof of Theorem \ref{thm.1}}
In this section, we prove Theorem \ref{thm.1} which consists of several steps.

\textbf{Step-1}.
Let us consider the energy functional $E \in C^1(X_0, \R)$ corresponding to the problem $(\mathcal{P})$ given by:
\begin{align} \label{Eq:I}
E(u)=\frac{1}{2}\norm{u}_{X_0}^2-\frac{\la}{2}|u|_2^2-\frac{\al}{q}|u|_q^q-\frac{\ba}{2^*}|u|_{2^*}^{2^*} , \, u \in X_0,
\end{align}
for $\la,\al \in \R, \ba>0$.

For $c,d \in \R, c \leq d,$ let us define

$$ E_c:=\Big\{u \in X_0:E(u) \geq c\Big\}, \,  E^d:=\Big\{u \in X_0: E(u) \leq d\Big\}. $$
Then note that $E^{-1}[c,d]= E_c \cap E^d$. 
For $t>0,$ let us define $B_t:=\{u\in X_0 :\|{u}\|_{X_0}<t \}.$
We now claim the following:

\textit{Claim 1}:
\textit{
 Let $\la>0.$ Then, there exists $l \in \N$,  $t, \mu >0$ such that:
	\begin{gather*}
	X_{l-1}^{\perp} \cap \pa B_t \subset E_\mu  \, \text{and}\,\,
	X_{l-1}^{\perp} \cap  (B_{t}\setminus\{0\}) \subset (E_0 \setminus E^{-1}(0)) .
	\end{gather*}
}
\textit{Proof of Claim 1}. By the definition of $\la^+$, we note that, there exists $l \in \N$ such that $\la^+=\la_l$. 
Using the Raleigh quotient characterization of $\la_l$ , Sobolev embedding and interpolation inequality, we have,
\begin{gather*}
	E(u)\geq \frac{1}{2}(1-\frac{\la}{\la^+})\norm{u}_{X_0}^2-c_1\norm{u}_{X_0}^p-c_2\norm{u}_{X_0}^{2^*}, \, \text{for all}\, u \in X_{l-1}^{\perp},
\end{gather*}
for some $c_1,c_2>0$. Now, the result follows from elementary analysis of the function
$$f(x)=\frac{1}{2}(1-\frac{\la}{\la^+})x^2-c_1x^p-c_2x^{2^*}, x \geq 0. $$
	Therefore, for $0<\|{u}\|_{X_0}<r, E(u) \geq f(\norm{u}_{X_0})>0 $ and we have the claim.
%
%
%
%



 \textbf{Step-2}. Let $\mu, l$ be given as in Step-1 and $\ba>0$. Let us denote 
 $$c:=\sup_{h \in \Ga} \inf_{u \in \pa(B_1) \cap X^{\perp}_{l-1}} E(h(u)),\quad b:=\inf_{k \in \Ga^*} \sup_{u \in K} E(u)
 \quad\mbox{and}\quad c^*:=\frac{s}{N}\frac{1}{\ba^{\frac{N-2s}{2s}}}S^{N/2s},$$
where 
$\Ga:=\Big\{h \in C(X_0,X_0):h(B_1) \subset
E_0 \cup \bar{B_r}, h \,\ \text{is an odd homeomorphism of}\, X_0 \Big\} $,
and $$\Ga^*:=\Big\{K \subset X_0 :K \,\ \text{is compact, symmetric and} \,\ \ga(K \cap h(\pa(B_1)) \geq l \,\text{for all}\,
h \in \Ga\Big \},$$ where $\ga$ is the genus. We notice that $\ga(K \cap h(\pa(B_1))$ in $\Ga^*$ is well-defined.

%



\textit{Claim 2}.
	\textit{
	The following holds true:
	$$0<\mu \leq c \leq b<c^*.$$ }
\textit{Proof of Claim 2}. As $h=r \, Id \in \Ga$, using Step-1, we have,
$$c \geq \inf_{{u \in \pa(B_1)} \cap X^{\perp}_{l-1}}E(ru) \geq \mu >0. $$
We observe that for any $K \in \Ga^*$ and $h \in \Ga, $ we have $K \cap h(\pa(B_1) \cap X^{\perp}_{l-1}) \neq \emptyset. $ See \cite{Amb-Rabin-73} for details.
Thus, for $w_{h,k} \in K \cap h(\pa B_1 \cap X^{\perp}_{l-1})$, there exists $u_{h,k} \in \pa(B_1) \cap X^{\perp}_{l-1}$ such 
	that $h(u_{h,k})=w_{h,k}. $ Therefore, 
	$$\inf_{u \in \pa B_1 \cap X^{\perp}_{l-1}}E(h(u)) \leq E(h(u_{h,k}))=E(w_{h,k}) \quad \mbox{for all}\quad h \in \Ga. $$
	
	This implies,
	$$c \leq \sup_{h \in \Ga}E(w_{h,k}) \leq \sup_{u \in K} E(u). $$ Hence, $c \leq b. $
	
	For any finite dimensional subspace $\widetilde{X} \subset X_0$ we have $E_0 \cap \widetilde X$ is bounded in $X_0. $
	Suppose not. Then,
	as in finite dimensional subspace $\widetilde X, $ all norms are equivalent, so there exists sequence $\{u_k\}_{k \geq 1} 
	\subset E_0 \cap \widetilde X$ such that $\norm{u_k}_{X_0} \to \infty, k \to \infty. $
	Now, as $u_k \in E_0, $ we have
	$$0 \leq E(u_k) \leq \frac{1}{2}\norm{u_k}_{X_0}^2-c_1\norm{u_k}_{X_0}^2-c_2\norm{u_k}_{X_0}^p-c_3\norm{u_k}_{X_0}^{2^*}.$$
	Letting $k \to \infty,$ we have the contradiction.
	In particular, we take $\widetilde X= V_\var. $ Then $E_0 \cap V_\var$ is bounded. So there exists $R>r>0$ sufficiently large such
	that $E_0 \cap V_\var \subset K_R$ where $K_R:=V_\var \cap \overline{B_R}. $ We note that $K_R$ is compact and symmetric. 
	Then for any $h \in \Ga, $ we have,
	$$h(B_1) \cap V_\var \subset (E_0 \cup \overline{B_R})\cap V_\var= (E_0 \cap V_\var) \cup K_R \subset  K_R.  $$
	Since $h(B_1) \cap V_\var$ is a bounded neighbourhood of $0$ in the $l$ dimensional
	subspace $V_\var,$ then $\ga(\pa(h(B_1) \cap V_\var))=l. $ As $h$ is a homeomorphism, we have 
	$$l \leq \ga(\pa(h(B_1) \cap V_\var)) \leq \ga(h(\pa B_1) \cap V_\var). $$
	Also $h$ being a homeomorphism, $h(B_1) \cap V_\var \subset K_R. $ This implies, $h(\pa B_1) \cap V_\var \subset K_R $ and $K_R \cap h(\pa B_1)=h(\pa B_1) \cap V_\var. $
	Therefore, $\ga(K_R \cap h(\pa B_1)) \geq l. $ So, $K_R \in \Ga^*.$ Thus 
	$$b \leq \sup_{u \in V_\var} E(u) \leq \sup_{u \in V_\var} E_{0 \ba}(u),$$ where $E_{0 \ba}(u)=\frac{1}{2}\norm{u}_{X_0}^2-\frac{\la}{2}|u|_{2}^2
	-\frac{\ba}{2^*}|u|_{{2^*}}^{2^*}. $
	 Then we have,
	\begin{gather*}
	E_{0 \ba}(u) \leq \max_t E_{0 \ba}(tu) \leq \frac{s}{N}\frac{1}{\ba^{\frac{N-2s}{2s}}} \bigg(\frac{\norm{u}_{X_0}^2-\la|u|^2_{2}}{|u|^2_{{2^*}}}\bigg)^{N/2s},
	\end{gather*}
	$u \in X_0 \setminus \{0\}$.
Indeed, by a rudimentary analysis follows by constructing the function
 $N(t)=\frac{t^2}{2}\norm{u}_{X_0}^2-\frac{\la t^2}{2}|u|_{2}^2-\frac{t^{2^*}\ba}{2^*}|u|_{{2^*}}^{2^*}, t \geq 0. $
So,	 $N'(t_0)=0$ implies $t_0=0$ or $$ t_0=\Bigg(\frac{\norm{u}_{X_0}^2-\la|u|^2_{2}}{\ba |u|^{2^*}_{{2^*}}}\Bigg)^\f{1}{2^*-2}. $$
	Hence, $N''(t_0)=-(2^*-2)\bigg(\norm{u}_{X_0}^2-\la|u|_{2}^2 \bigg) \leq 0. $
	So, $N(t)$ attains maximum at $t=t_0. $ Therefore, an easy computation yields,
	\Bea
	N(t_0) &=& \frac{t_0^2}{2}(\norm{u}_{X_0}^2-\la|u|_{2}^2-t_0^{(2^*-2)}\frac{2\ba}{2^*}|u|_{{2^*}}^{2^*})\\
	&=& \frac{s}{N}\frac{1}{\ba^{\frac{2}{2^*-2}}}\bigg(\frac{\norm{u}_{X_0}^2-\la|u|^2_{2}}{|u|^{2}_{{2^*}}}\bigg)^{\f{2^*}{2^*-2}}.\\
	\Eea
	This implies,
	$$E_{0 \ba}(u) \leq \max_t E_{0 \ba}(tu) \leq \frac{s}{N}\frac{1}{\ba^{\frac{2}{2^*-2}}}\bigg(\frac{\norm{u}_0^2-\la|u|^2_{L^2(\Om)}}{|u|^{2}_{L^{2^*}(\Om)}}\bigg)^{\f{N}{2s}}.$$
	
	So by (\ref{Eq:Q-la})  we have, 
	$$\sup_{u \in V_\var}E_{0 \ba}(u) \leq \frac{s}{N}\frac{1}{\ba^{\frac{N-2s}{2s}}}S^{N/2s}=c^*. $$ This finishes the proof of Claim 2.
	
	\textbf{Step-3}. In this step, we show $E$ satisfies $(PS)_c$ condition for any $c <c^*. $
	 To this goal, let $\{u_k\}_{k \geq 1}$ be a sequence in $X_0(\Om)$ such that $E(u_k) \to c$ and $E'(u_k) \to 0$ in $(X_0(\Om))',
	k \to \infty. $ It is easy to check that $\{u_k\}_{k \geq 1}$ is bounded in $X_0(\Om)$.  Hence, 
	going to a subsequence, if necessary, we can assume that as $k \to \infty, $ there exists $u \in X_0(\Om)$ such that
	\begin{align}\label{*1}
		&u_k \deb u\quad\mbox{in}\quad X_0(\Om), \notag\\
		&u_k \To u\quad\mbox{strongly in}\quad L^q(\Rn)\quad\mbox{for}\quad q \in [1,2^*),\notag\\
		&u_k \To u\quad\mbox{a.e. in}\quad\Rn\quad\mbox{for}\quad 1\leq \ga < 2^*,
	\end{align}
	and there exists$\quad l \in L^q(\Rn)\quad\mbox{such that}$
	\begin{align*}
		&|u_k(x)|\leq l(x)\quad\mbox{a.e. in}\quad \Rn\quad\mbox{for all}\quad k \geq1.
	\end{align*}
	
	Using Vitali Convergence Theorem, we have,
	\begin{gather*}
	\Iom u_k \phi \dx \to \Iom u \phi \dx, k \to \infty \quad\mbox{for all}\quad \phi \in X_0(\Om), 
	\end{gather*}
	\begin{gather*}
	\Iom |u_k|^{p-2}u_k \phi \dx \to \Iom |u|^{p-2}u \phi \dx,k \to \infty \quad\mbox{for all}\quad \phi \in X_0(\Om),
	\end{gather*}
	 and
	\begin{gather*}
	\Iom |u_k|^{2^*-2}u_k \phi \dx \to \Iom |u|^{2^*-2}u \phi \dx,k \to \infty \quad\mbox{for all}\quad \phi \in X_0(\Om). 
	\end{gather*}
	Hence, $u$ is the solution of :
	$\left(-\Delta\right)^s u-\la u = \al |u|^{p-2}u +  \ba|u|^{2^*-2}u$. 
	This gives,
	$$\|{u}\|_{X_0}^2-\la|u|_{2}^2-\al|u|_{p}^p-\ba|u|^{2^*}_{{2^*}}=0. $$
	Using this, we have
	\begin{align} \label{*2}
		E(u)=\al(\frac{1}{2}-\frac{1}{p})|u|_{p}^p +\ba(\frac{1}{2}-\frac{1}{2^*})|u|_{{2^*}}^{2^*}
	\end{align}
	Let $v_n=u_n-u. $ Brezis-Lieb lemma leads to:
	\begin{align} \label{*3}
		|u_n|_{{2^*}}^{2^*}=|u|_{{2^*}}^{2^*}+|v_n|_{{2^*}}^{2^*}+o(1).
	\end{align}
	So we have,
	\begin{equation*}
	\begin{aligned}
	&E(u)+\frac{1}{2}\|{v_n}\|_{X_0}^2-\frac{\ba}{2^*}|v_n|^{2^*}_{{2^*}} \\
	&= \frac{1}{2}(\|{u}\|_{X_0}^2-\|{u_n}\|_{X_0}^2)
	-\frac{\la}{2}(|u|_{2}^2-|u_n|^2_{2})-\frac{\al}{p}(|u|^p_{p}-|u_n|^p_{p})\\
	 &-\frac{\ba}{2^*}(|u|^{2^*}_{{2^*}}-|u_n|^{2^*}_{{2^*}}+|v_n|^{2^*}_{{2^*}})+\frac{1}{2}(\|{u_n}\|_{X_0}^2-\|{u}\|_{X_0}^2)-\<u_n,u\>+E(u_n)\\
	&= \|{u}\|_{X_0}^2-\<u_n,u\>-\frac{\al}{p}(|u|^p_{p}-|u_n|^p_{p})-\frac{\la}{2}(|u|_{2}^2-|u_n|^2_{2})+o(1)+E(u_n)\\
	&= o(1)+E(u_n) \to c, n\to \infty.
	\end{aligned}
	\end{equation*}
	This implies,
	\bea \label{*4}
	E(u)+\frac{1}{2}\|{v_n}\|_{X_0}^2-\frac{\ba}{2^*}|v_n|^{2^*}_{{2^*}} &\to& c, n\to \infty.
	\eea
	Since $\<E'(u_n),u_n\> \to 0,n \to \infty$ so we have,
	\begin{equation*}
	\begin{aligned}
	\|{v_n}\|_{X_0}^2-\ba|v_n|_{{2^*}}^{2^*} &= \|{u_n}\|_{X_0}^2-\la|u_n|^2_{2}-\al|u_n|^p_{p}-\ba|u_n|^{2^*}_{{2^*}}+\la|u_n|^2_{2}+\al|u_n|_{p}^p
	\\
	 &+\ba(|u_n|^{2^*}_{{2^*}}-|v_n|^{2^*}_{{2^*}})+\|{u}\|^2-2\<u_n,u\>\\
	&= \<E'(u_n),u_n\>+\la|u_n|^2_{2}+\al|u_n|_{p}^p+\ba(|u|_{{2^*}}^{2^*}+o(1))+\|{u}\|_{X_0}^2-2\<u_n,u\>\\
	&\to \la|u|_{2}^2+\al|u|_{p}^p+\ba|u|_{{2^*}}^{2^*}-\|{u}\|_{X_0}^2\\
	&=-\<E'(u),u\>= 0.
	\end{aligned}
	\end{equation*}
This gives us,
	\bea \label{*5}
	\|{v_n}\|_{X_0}^2-\ba|v_n|_{{2^*}}^{2^*} \to 0.
	\eea
	Since $(u_n)_{n \geq 1}$ is bounded in $X_0(\Om), $ so $\|{v_n}\|_{X_0}$ is bounded. Hence, upto a subsequence we may assume that
	$\|{v_n}\|_{X_0}^2 \to b. $
	Using $(\ref{*5})$ we have,
	$$\ba|v_n|^{2^*}_{{2^*}} \to b,$$ and by Sobolev embedding,
	$$\|{v_n}\|_{X_0}^2 \geq S|v_n|^{2^*}_{{2^*}}.$$These two estimates imply
	\bea \label{*6}
	b \geq S(\frac{b}{\ba})^{2/{2^*}}.
	\eea
	If $b=0, $ then we are done. 
	Suppose not. Then we have 
	$b \geq \frac{S^{\frac{2^*}{2^*-2}}}{\ba^{\frac{2}{2^*-2}}}. $
	From $(\ref{*3})$ we have,
	$E(u)+\frac{1}{2}b-\frac{1}{2^*}b=c. $ 
	This implies,
	\Bea
	c &=& (\frac{1}{2}-\frac{1}{2^*})b+E(u)\\
	&=& \frac{s}{N}(b+\ba|u|_{{2^*}}^{2^*})+(\frac{1}{2}-\frac{1}{p})|u|_{p}^p\\
	&\geq&\frac{s}{N}\frac{S^{\frac{2^*}{2^*-2}}}{\ba^{\frac{2}{2^*-2}}}+\frac{s}{N}\ba|u|^{2^*}_{{2^*}}+\al(\frac{1}{2}-\frac{1}{p})|u|_{p}^p\\
	&\geq& c^*+\frac{s}{N}\ba|u|^{2^*}_{{2^*}}+\al(\frac{1}{2}-\frac{1}{p})|u|_{p}^p.
	\Eea
	As $p \geq 2, \al, \ba >0, $ we have $c \geq c^*, $ which is a contradiction to Claim 2 in Step-2.
	Hence, $b=0$ and the conclusion follows. This finishes Step-3.

\textbf{Step-4}. In this Step, we conclude the proof of Theorem \ref{thm.1}.  Before that, we  recall the following lemma from [Lemma 3.1, \cite{Willem-Minimax}]. (Also, see [Lemma 2.3, \cite{Yue-Zou}]).
\begin{lemma} \label{lemma3}
	Let $\mu, c$ as in Step-1. Let $\var \in (0,\frac{1}{2}(c-\frac{\mu}{2})), \de>0.$ Let $h \in \Ga$ be such that 
	$$\inf_{u \in \pa B_1 \cap E^{\perp}_{l-1}}E(h(u)) \geq c-\var. $$
	Then, there exists $v_\var \in X_0$ such that
	\begin{enumerate}
		\item 
		$c-2\var \leq E(v_\var) \leq c+2\var,$
		\item
		$dist(v_\var,h(\pa B_1 \cap Y^{\perp}_{l-1}))\leq 2\de,$
		\item
		$\norm{E'(v_\var)} \leq \frac{8 \var}{\de}.$
	\end{enumerate}
\end{lemma}

Using Lemma \ref{lemma3}, we can say there exists a sequence $(v_n)_{n \geq 1} \in X_0$ such that $E(v_n) \to c$ and
$E'(v_n) \to 0$, as $n\to \infty$. So by Step-3,  there exists a subsequence, still denoted by $v_n$ and $u \in X_0$
such that $v_n \to u $ in $X_0$ and $E(u)=c, \, E'(u)=0, $ that is, $\<E'(u),\phi\>=0$ for all $\phi \in X_0$.
In particular, for $\phi=\phi_1>0$ where $\phi_1$ is the first eigen function of $\left(-\Delta\right)^s, $
we have, $$\<E'(u),\phi_1\>=0. $$ This implies,
$$\<u,\phi_1\>_{X_0}-\la\Iom u\phi_1=\al \Iom|u|^{p-2}u \phi_1 + \ba \Iom |u|^{2^*}u \phi_1. $$
Note that $\<u,\phi_1\>_{X_0}=\<u,\left(-\Delta\right)^s\phi_1\>_2=\la_1\<u,\phi_1\>_2=\la_1\Iom u \phi_1. $
These two estimates together imply,
\bea \label{la_1}
(\la_1-\la)\Iom u\phi_1=\al \Iom |u|^{p-2}u \phi_1 + \ba \Iom |u|^{2^*}u \phi_1.
\eea
Let $\la_1 \leq \la, \al, \ba >0. $ Then $u^+ \neq 0, u^- \neq 0. $
If not, that is, if $u^+\equiv 0, $ then R.H.S of \eqref{la_1} $>0$ but L.H.S. $<0. $ and if  
$u^-\equiv 0, $ then R.H.S of \eqref{la_1} $<0$ while L.H.S. $>0. $ Hence, $u$ is sign changing and this completes the proof.

\section{Proof of Theorem \ref{thm.2}}
In this section, we consider the case $N \geq 3s$, $s \in (0,1)$, $\la \in \R, \al,\ba>0$. We rewrite the energy functional taking care of $\ba$ as a perturbation term,
$$E_{\al\ba}(u) \equiv E(u)=\bigg\{\frac{1}{2}\norm{u}_{X_0}^2-\frac{\la}{2}|u|_2^2-\frac{\al}{p}|u|_p^p\bigg\}
+\bigg\{-\frac{\ba}{2^*}|u|_{2^*}^{2^*}\bigg\}\\
:\equiv K_\al(u)+J_\ba(u),$$
where $K_\al(u)=\frac{1}{2}\norm{u}_{X_0}^2-\frac{\la}{2}|u|_2^2-\frac{\al}{p}|u|_p^p,\,J_\ba(u)=-\frac{\ba}{2^*}|u|_{2^*}^{2^*}. $\\
Let $A \subset X_0(\Om)$ be a subset. We denote the neighbourhood of $A$ by $A^d$ where 
$$A^d:=\bigcup_{u \in A}B_d(u),\,\,B_d(u):=\Big\{v \in X_0(\Om): \|{u-v}\|_{X_0} \leq d\Big\}. $$

The proof is also divided into several steps.

\textbf{Step-1}. In the first step,
we recall the following results from [Theorem 3.7, \cite{Willem-Minimax}]. 
 \begin{lemma} \label{lem-1}
	For each $k \in \N, \quad\text{there exists}\quad R_k>0$ such that $K_\al(u)<0$ with $u \in X_k, \|u\|_{X_0}=R_k. $
\end{lemma}
\begin{lemma} \label{lem-2}
	For each $k \in \N, $ there exists $0<r_k<R_k$ such that for $u \in X^{\perp}_{k-1}$ and $\norm{u}_{X_0}=r_k, K_\al(u) \to \infty,
	k \to \infty. $
\end{lemma}

 We also list some notations which we will use in this section:
$$B_k:=\Big\{u \in X_k:\norm{u}_{X_0} \leq R_k\Big\},
\Ga_k:=\Big\{h \in C(B_k,X_0(\Om)):h \, \text{is odd}\, , h|_{\pa B_k}=id\Big\},$$
$$c_{\al k}:=\inf_{h \in \Ga_k} \sup_{u \in B_k} K_\al(h(u)),
a_{\al k}:=\sup_{u \in \pa B_k}K_{\al}(u), b_{\al k}:=\inf_{u \in \pa Z_k} K_{\al}(u),$$ where
$$Z_k:=\Big\{u \in X^{\perp}_{k-1}, \|{u}\|_{X_0} \leq r_k\Big\}.
$$ \\
By Lemma \ref{lem-1}, there exists $k_1$ such that for all $k \geq k_1, $ there exists $0<r_k<R_k$ such that 
$u \in X^{\perp}_{k-1}, \|{u}\|_{X_0}=r_k, K_\al(u) \geq 1. $
Hence, for all $k \geq k_1, $ we have
$$\inf_{u \in X^{\perp}_{k-1} \|{u}\|_{X_0}=r_k}K_\al(u) \geq 1>0>K_\al(u) \quad\text{for all}\quad u \in X_k \,\text{with}\, \|{u}\|_{X_0}=R_k. $$
So, we have for $k$ large enough $b_{\al k} > a_{\al k}$ and $c_{\al k} \geq b_{\al k}. $
We refer [Theorem 3.5, \cite{Willem-Minimax}] for proof. 
For all $k \in \N$ large enough the functional $K_\al$ satisfies $(PS)_{c_{\al k}}$ condition.
We know that for any $c>0, K_\al$ satisfies $(PS)_c$ condition and thus, $K_\al$ has infinetly many critical values. For each
$k \in \N,$ let us define:
$$V_{\al k}:=\Big\{u \in X_0(\Om)\setminus \{0\}:K'_\al(u)=0,K_\al(u)=c_{\al k}\Big\}. $$
Therefore, $V_{\al k}$ is non-empty and compact. Let $S_{1k}:=\sup_{u \in V_{\al k}}\|{u}\|_{X_0}. $

\textbf{Step-2}. In this step, we prove the following claim.\\
\textit{Claim 1}.
\textit{
	For $c>0$, the following holds true:
	\begin{equation*}
	\lim_{\ba \to 0} \sup_{\|{u}\|_{X_0} \leq c}|J_\ba(u)|=\lim_{\ba \to 0} \sup_{\|{u}\|_{X_0} \leq c}|J'_\ba(u)|=0. 
	\end{equation*}
}
\textit{Proof of Claim 1}.
	We have,
	$$|J_\ba(u)|=\frac{\ba}{2}|u|^{2^*}_{2^*} \leq \frac{\ba}{2}S^{-2^*/2}\|{u}\|_{X_0}^{2^*}. $$
	This implies,
	$$\sup_{\|{u}\|_{X_0} \leq c} |J_\ba(u)| \leq \frac{\ba}{2}S^{-2^*/2}c^{2^*}.$$ and so
	$\lim_{\ba \to 0} \sup_{\|{u}\|_{X_0} \leq c}|J_\ba(u)|=0. $
	Let $v \in X_0(\Om). $ Then, we have,
	\Bea
	|\<J'_\ba(u),v\>| 
	&\leq& \ba \|{u}\|_{X_0}^{2^*-1}\|{v}\|_{X_0}S^{-2^*/2}.
	\Eea
	This evidently implies,
	\begin{gather*}
	\frac{\big|\<J'_\ba(u),v\>\big|}{\|{v}\|_{X_0}} \leq \ba \|{u}\|_{X_0}^{2^*-1}S^{-2^*/2}. 
	\end{gather*}
	Therefore, $\sup_{\|{u}\|_{X_0} \leq c}\|{J'_{\ba}(u)}\|_{(X_0)'} \leq \ba c^{2^*-1}S^{-2^*/2}$ and Claim 1 follows. 

\textbf{Step-3}. In this section, we again set some notaions:
 $$c_{\ba k}:=\inf_{u \in \Ga_k} \sup_{u \in B_k}E_{\al\ba}(h(u)), 
a_{\ba k}:=\sup_{u \in \pa B_k}E_{\al\ba}(u),  b_{\ba k}:=\inf_{u \in \pa Z_k} E_{\al\ba}(u). $$
Note that $b_{\ba k} >a_{\ba k}$ for $\ba$ sufficiently small and $k \in \N$ large enough. 
Indeed, $E_{\al\ba}(u)=K_\al(u)+J_\ba(u) \leq K_{\al}(u)$ for all $u \in X_0(\Om). $
So, $a_{\ba k} \leq a_{\al k}<0. $ For $\ba >0$ small enough, we have
\Bea
b_{\ba k}&=& \inf_{u \in \pa Z_k}(K
_\al(u)+J_\ba(u))
> a_{\ba_k}.
\Eea
Also, for any $h \in \Ga_k$ we have, $h(B_k) \cap \pa Z_k \neq \emptyset. $ So, $c_{\ba k} \geq b_{\ba_k}$ by 
[Theorem 3.5,\cite{Willem-Minimax}] and we get a $(PS)_{c_{\ba k}}$ sequence of $I_{\al \ba}$. 

\textbf{Step-4}. In this section, we claim the following. 

\textit{Claim 2}.
\textit{
	For $\la \in \R, \ba,\al >0, \{u_{\ba kn}\}_{n \geq 1}$ is a $(PS)_{c_{\ba k}}$ sequence of $E_{\al \ba}, $ then there exists
	$S_{2k} >0$ independent of $\ba$ such that $\|{u_{\ba kn}}\|_{X_0} \leq S_{2k}$ for all $n$ large enough. }
	
\textit{Proof of Claim 2}.
	We note that,
	$c_{\ba k} \leq c_{\al k}$, since  $E_{\al \ba}(u) \geq K_\al(u)$ for all $u \in X_0(\Om)$.
	By choosing $\rho \in (\frac{1}{p},\frac{1}{2})$ for $n$ large enough we have,
	\Bea
	c_{\al k}+1+\|{u_{\ba kn}}\|_{X_0} &\geq& c_{\ba k}+1+\|{u_{\ba kn}}\|_{X_0} \\
	&\geq& E_{\al \ba}(u_{\ba kn})-\rho \<E'_{\al \ba}(u_{\ba kn}),u_{\ba kn}\> \\
	&\geq&(\frac{1}{2}-\rho)\|{u_{\ba kn}}\|_{X_0}^2+\al((\rho-\frac{1}{p})-\var)|u_{\ba kn}|_p^p
	-c_\var,
	\Eea
	where $0<\var<\al(\rho-\frac{1}{p})$ and $c_\var$ independent of $\ba. $ This completes the proof of the Claim.

\textbf{Step-5}. We signify $\tilde{S_k}:=2\max\{S_{1k},S_{2k},R_k\}. $
In this step, we prove the following Claim. 

\textit{Claim 3}.
\textit{
	We have $lim_{\ba \to 0}c_{\ba k}=c_{\al k}. $}
	
\textit{Proof of Claim 3}.
	Let $\var>0. $ From the definition of $c_{\al k}$ we have $h_\var \in \Ga_k$ such that $$\sup_{v \in B_k} K_\al(h_\var(v)) 
	\leq c_{\al k} + \eps, $$ that is,
	$$\sup_{u \in h_\var(B_k)}K_\al(u)-\var \leq c_{\al k}. $$
	Let $\Ga_{k,\var_j \al}:=\{h_\var \in \Ga_k:\sup_{u \in h_\var (B_k)}K_\al(u)-\var \leq c_{\al k}\}. $ Then, 
	$$\lim_{\var \to 0} \inf_{h \in \Ga_{\var,k,\al}} \sup_{u \in h(B_k)}K_\al(u)=c_{\al k}. $$
	
	We notice that, for each $\var>0, $ there exists $v_\var \in h_\var(B_k)$ such that $$\sup_{u \in h_\var(B_k)}K_\al(u)=K_\al(v_\var),$$
	and $\{v_\var\}_{\var \geq 0}$ is a $(PS)_{c_{\al k}}$ sequence of $K_\al$ as $\var \to 0. $
	

Therefore, for $\var>0$ small enough we have $v_\var \in B_{\widetilde{S_k}}(0). $ Thus, 
$$c_{\al k} \leq K_\al(v_\var) \leq \sup_{u \in h_\var(B_k) \cap B_{\widetilde {S_k}}(0)} K_\al (u). $$
So, we have, $c_{\al k}=\lim_{\var \to 0} \inf_{h \in \Ga_{k,\var,\al}} \sup_{u \in h(B_k) \cap B_{\widetilde {S_k}}(0)} K_\al (u). $
Similarly, we can show that 
$$c_{\ba k}=\lim_{\var \to 0} \inf_{h \in \Ga_{k,\var,\ba}}\sup_{u \in h(B_k) \cap B_{\widetilde{S_k}}(0)}J_\ba(u),$$
where $$\Ga_{k,\var,\ba}:=\{h \in \Ga_k:\sup_{u \in h(B_k)}E_{\al \ba}(u)-\var \leq c_{\ba k}\}. $$

By Claim 2, for any $h \in \Ga_k$ and $u \in h(B_k) \cap B_{\widetilde{S_k}}(0)$ we get,
$K_\al(u)=\lim_{\ba \to 0} E_{\al \ba}(u). $ As $c_{\al k} \geq c_{\ba k}$ for any $h \in \Ga_{k,\var,\ba}$ such that
$$\sup_{u \in h(B_k)\cap B_{\widetilde{S_k}}(0)} E_{\al \ba}(u)-\var \leq c_{\ba k} \leq c_{\al k}. $$
So we have $h \in \Ga_{k,\var,\al}$ as $\ba \to 0. $ Hence,
$\lim_{\ba \to 0} \inf c_{\ba k} \geq c_{\al k}. $ Therefore, $c_{\al k}=\lim_{\ba \to 0} c_{\ba k}. $ This completes the proof of Claim-3.

%


\textbf{Step-6}.
In this step, we claim the following, thereby completing the proof. 

\textit{Claim 4}.
\textit{For any $d>0, $ there exists $\ba_{0k}>0$ such that for any $\ba \in (0,\ba_{0k})$ and any $(PS)_{c\ba k}$ sequence $\{u_{\ba k n\}_{n \geq 1}}$
of $E_{\al \ba}, $ there exists $n_0 \in \N$ such that $u_{\ba kn} \in V_{\al k}^d$ for all $n \geq n_0. $}

\textit{Proof of Claim 4}. We will prove the claim by the method of contradiction. Suppose not, then there exists $d >0, $ any sequence
$\ba_j>0$ as $n_j \to \infty, $ such that 
\begin{equation*}
\lim_{j \to \infty}E_{\al \ba_j}(u_{\ba_j k_{n_j}})=c_{\ba_j k}, \, \lim_{j \to \infty}E'_{\al \ba_j}(u_{\ba_j k_{n_j}})=c_{\ba_j k},
\end{equation*}
but $u_{\ba_j k_{n_j}} \notin V_{\al k}. $
Using Claim 2, there exists $N_{2k}>0$ such that $\|{u_{\ba j k_{n_j}}}\|_{X_0} \leq N_{2k}. $
Using Claim 1 and Claim 3, we have $\lim_{j \to \infty}K'_{\al} (u_{\ba j k_{n_j}})=c_{\al k}. $
Indeed, 
\Bea
|K_\al (u_{\ba j k_{n_j}})-c_{\al k}|&\leq& |E_{\al \ba_j}(u_{\ba j k_{n_j}})-c_{\ba_jk}|+|c_{\ba_j k}-c_{\al k}|+|J_{\ba_j}(u_{\ba j k_{n_j}})|\\
&\to& 0, j \to \infty,
\Eea
and
\Bea
\|{K'_\al (u_{\ba j k_{n_j}})}\| &\leq&  \|{K'_\al (u_{\ba j k_{n_j}})+J'_{\ba_j} (u_{\ba j k_{n_j}})}\|+\|{J'_{\ba_j} (u_{\ba j k_{n_j}})}\|\\
&=& \|{E'_{\al \ba}(u_{\ba j k_{n_j}})}\|+\|{J'_{\ba_j}(u_{\ba j k_{n_j}})}\| \\
&\to& 0, j \to \infty.
\Eea
Therefore, $\{u_{\ba j k_{n_j}}\}$ is a $(PS)_{c_{\al k}}$ sequence of $K_\al. $ As $K_\al$ satisfies $(PS)_{c_{\al k}}$ condition, so up to a subsequence, there exists $u_{0k} \in V_{\al k}$ such that $u_{\ba j k_{n_j}} \to u_{0k}. $
Therefore, $u_{\ba j k_{n_j}} \in B_d(u_{0k}) \subset V^d_{\al k}$ for j large, which is a contradiction. This proves Claim 4. 

By the above Claim 4, we consider $d_k>0$ small so that for any $u \in V_{\al k}^{2d_k}, $ we have $u \nequiv 0. $
For any $\ba \in (0,\ba_{0k}), $ let $\{u_{\ba k_n}\}$ be $(PS)_{c \ba_k}$ sequence for $E_{\al \ba}, $ then $\{u_{\ba k_n}\} \in V^{d_k}_{\al_k}. $ 
Therefore, there exists $u_{0kn} \in V_{\al k}$ such that $u_{\ba k_n} \in B_{d_k}(u_{0k_n}). $
As $V_{\al k}$ is compact, there exists $u_{0k} \in V_{\al k}$ such that $u_{0k_n} \to u_{0k}. $ So, we have $u_{\ba k_n} \in B_{2d_k}(u_{0k})$
for large $n. $ Also, there exists $u_{\ba k}$ such that $u_{\ba k_n} \deb u_{\ba k}$ and $E'_{\al \ba}(u_{\ba k})=0. $
As $\ba_{2d_k}(u_{0k})$ is closed and convex, so it is weakly closed. Therefore, $u_{\ba k} \in B_{2d_k}(u_{0k}). $ Hence, $u_{\ba k}$ is a nontrivial
critical point of $E_{\al \ba} $ and we have for any $m \in \N, $ for every $\ba \in (0,\ba_m) $ the problem $(\mathcal{P})$ has $m$ nontrivial solutions
where $\ba_m:=\min_{1 \leq k \leq n} \ba_{0k}. $ This completes the proof. 

\section{Proof of Theorem \ref{thm.3}}
Let us denote the energy functional as in Step-1 of the proof of Theorem \ref{thm.1},
$$
E_{\al}(u)=E(u)=\frac{1}{2}\|u\|_{X_0}^2-\frac{\la}{2}|u|^2_2-\frac{\ba}{2^*}|u|_{2^*}^{2^*}
-\frac{\al}{q}|u|_q^q:=E_0(u)+G_\al(u),
$$ where
 $$E_0(u)=\frac{1}{2}\|u\|_{X_0}^2-\frac{\la}{2}|u|_2^2-\frac{\ba}{2^*}|u|_{2^*}^{2^*}\,
\text{and}\, G_{\al}(u)=-\frac{\al}{q}|u|_q^q ,$$ with
 $\la>0, \al<0, \ba>0, N \geq 4s$ with $s \in (0,1)$.

We prove the result in the following steps. \\
\textbf{Step-1}. Approaching similarly as in the proof of Claim 1 of Step-1 in the proof of Theorem \ref{thm.1}, using Raleigh quotient characterization of $\la_l$,
Sobolev embedding that,
$$
E_0(u) \geq \frac{1}{2} \bigg( 1-\frac{\la}{\la^+}\bigg) \|u||_{X_0}^2-C\|u\|_{X_0}^2,   
$$
for all $u \in X_{l-1}^{\perp}$ and for some $C>0$.
This yields that, there exists $\mu_0, r_0>0$ such that,
$$
E_0(u) >\mu_0 \,\text{for all}\, u \, \text{with}\, \|u\|_{X_0}=r_0,
$$
and
$$
E_0(u) > 0 \,\text{for all}\, u \, \text{with}\, 0<\|u\|_{X_0}<r_0.
$$
Let $\Ga_0=\bigg\{ h \in C(X_0,X_0) : h(B_1) \subset (E_0)_0 \cup \overline{B}_{r_0},\, h \,\text{is an odd homeomorphism of}\,  X_0     \bigg\}$, where $E_0=\big\{ u \in X_0: E_0(u) \geq 0 \big\}$ and $B_{r_0}=\big\{ u \in X_0: \|u\|_{X_0}<r_0   \big\}$, and
$$
\Ga_0^*=\bigg\{ K \subset X_0: K \, \text{is compact, symmetric and}\, \ga(K \cap h(\pa B_1)) \geq l \,\text{for all}\, h \in \Ga \}, 
$$
where $\Ga$ is the genus. We set the following notations:
$$
b_0:=\inf_{K \in \Ga_0^*} \sup_{u \in K} E_0(u), \, c_0:=\sup_{h \in \Ga_0} \inf_{u \in \pa B_1 \cap X_{l-1}^{\perp}} E_0(h(u)).
$$
With tha above notations, proceeding as in Claim 2 of Step-2 in the proof of Theorem \ref{thm.1} that,
$$ 0<\mu_0 \leq c_0\leq b_0<c^*, $$
and similarly by Step-3 in the proof of Theorem \ref{thm.1}, we get,
$E_0$ has a (PS)$_{c_0}$ sequence and $c_0$ is a critical value of $E_0$.
Let us consider the set of critical points of $F$ with critical value $c_0$, that is, 
$$
Z_{E_0}:=\bigg\{u \in X_0 \setminus \{0\} : E_0(u)=c_0,\, E_0'(u)=0  \bigg\}.
$$
Hence, the set $Z_{E_0} \neq \emptyset$ and compact. We now signify:
$L_1:=\sup_{u \in Z_{E_0}}\|u\|_{X_0}$.    

\textbf{Step-2}. As in Claim 2 of Step-2 in the proof of Theorem \ref{thm.2}, we can show that for any $C>0$,
$$
\lim_{\al \to 0} \sup_{\|u\|_{X_0} \leq C} |G_{\al}(u)|=0=\lim_{\al \to 0} \sup_{\|u\|_{X_0} \leq C}  |G'_{\al}(u)|.
$$
Let us set
$\Ga_\al=\bigg\{ h \in C(X_0,X_0) : h(B_1) \subset (E_\al)_0 \cup \overline{B}_{r_0},\, h \,\text{is an odd homeomorphism of}\,  X_0  \bigg\}$.
 We notify:
$c_\al=\sup_{h \in \Ga_\al} \inf_{u \in \pa B_1 \cap X_{l-1}^{\perp}} E_\al (h(u))$. Then, one can observe that,
$$ E_\al (u) \geq E_0(u) \geq \mu_0 \,\text{for} \, u \in X_{l-1}^{\perp}, \|u\|_{X_0}=r_0.
$$
Hence, $c_\al \geq \mu_0$ and $E_\al$ satisfies (PS)$_{c_\al}$.

\textbf{Step-3}. In this step, we first claim the following.

\textit{Claim 1}. \textit{Let $\la>0, \al<0$ with $|\al|$ small, $\ba>0$. Assume $\{(u_\al)n \}_n$ be a (PS)$_{c_\al}$ sequence of $E_\al$. Then, there exists $L_2>0$ (independent of $\al$) such that,
	$$  \| (u_\al)_n\|_{X_0} \leq L_2 \, \,\text{for all}\,\, n \in \N.   $$
}

\textit{Proof of Claim 1}. For $\al>0,$ let us denote the set
$$
\Ga_{\al}^*=\big\{ K \subset X_0: K \,\text{is compact, symmetric and}\, \ga(K \cap h(\pa (B_1))) \geq l \,\text{for all}\, h \in \Ga_{\al}  \big\}. 
$$
We set $b_{\al}=\inf_{K \in \Ga_\al^*} \sup_{u \in K} E_\al(u)$. We observe that, $c_\al \leq b_\al$. There exists $R_0>r_0>0$ large enough such that
for $|\al|$ small, we have,
$(E_\al)_0 \cap V_\eps \subset \overline{B}_{R_0}$. Let $K_{R_0}:=V_\eps \cap \overline{B}_{R_0}$. Then, $K_{R_0} \subset \Ga_{\al}^*$. 

Thus, we obtain,
$$
b_\al \leq \sup_{u \in K_{R_0}} E_{\al}(u) \leq \sup_{V_\eps} E_0(u)+c_0<c^*+c_0=c^*,
$$
where $c_0>0$ (independent of $\al$). By choosing $\rho \in (\frac{1}{2^*}-\frac{1}{q})$, as in the proof of Claim 2 of Step-2 in the proof of Theorem \ref{thm.2}, we note that,
\begin{equation}
\begin{aligned}
c^*+\|u_{\al_n}||_{X_0} &\geq E_{\al}(u_{\al n}) -\rho \< E'_{\al}(u_{\al n}), u_{\al_n}\>\\
&\geq \big(\frac{1}{2}-\rho \big) \|u_{\al_n}\|_{X_0}^2+\bigg( \ba \big( \rho-\frac{1}{2^*} \big)-\bar{\eps}  \bigg)|u_{\al_n}|_{2^*}^{2^*}-c_{\bar{\eps}}\\
&\geq \big(\frac{1}{2}-\rho \big) \|u_{\al_n}\|_{X_0}^2-c_{\bar{\eps}},
\end{aligned}
\end{equation}
where $\bar{\eps} \in (0,\ba(\rho-\frac{1}{q}))$, $c_{\bar{\eps}}$ is independent of $\al$. This yields, there exists $L_2>0$ such that 
$$
\|u_{\al_n}\|_{X_0} \leq L_2 \,\text{for all}\, n \in \N.
$$

\textbf{Step-4}. We note that, $\Ga_0 \subset \Ga_{\al}$ for $\al<0$. We also observe that, for $u \in X_0$, we have,
$$
\lim_{\al \to 0}E_{\al}(u)=E_0(u).
$$
Therefore, for any $h \in \Ga_{\al}$ with $E_\al(h(B_1)) \geq 0,$ we have,
$E_0(h(B_1)) \geq 0$. In a similar fashion, as in Claim 3 of Step-5 in the proof of Theorem \ref{thm.2}, we obtain, 
\begin{gather}\label{eq:lim}
\lim_{\al \to 0}c_\al=c_0.
\end{gather}
Using similar arguments as in Step-6 in the proof of Theorem \ref{thm.2}, for any $d>0$, there exists $\al_0>0$ such that for any $\al \in (-\al_0,0)$, there exists (PS)$_{c_\al}$ sequence of $E_\al$, denoted by 
$\{u_{\al_n}\} \subset Z_{E_0}^d$. Taking $d>0$ sufficiently small, we obtain a nontrivial solution $u_\al$ in a neighborhood of a solution $u_0$
of $E_0$ with $E_0(u_0)=c_0$. This finishes the proof.

\section{Acknowledgement}
The author wishes to thank Prof. Mousomi Bhakta for giving helpful suggestions and comments. This research is supported by the Czech Science Foundation, project GJ19--14413Y.

\bibliography{biblio}
\bibliographystyle{acm} 

\end{document}